\documentclass[12pt]{amsart}
\usepackage[colorlinks=true]{hyperref}
\usepackage{amsmath}
\usepackage[english,  activeacute]{babel}
\usepackage[utf8]{inputenc}
\usepackage{amssymb}
\usepackage{amsthm}
\usepackage{graphics,graphicx}
\usepackage{enumerate}
\usepackage{array}
\usepackage{bm}
\usepackage{a4wide}
\usepackage{color, url}
\usepackage{float}
\usepackage{multirow}

\theoremstyle{plain}
\newtheorem{theorem}{Theorem}[section]
\newtheorem{proposition}[theorem]{Proposition}
\newtheorem{corollary}[theorem]{Corollary}
\newtheorem{lemma}[theorem]{Lemma}

\theoremstyle{definition}
\newtheorem{definition}[theorem]{Definition}
\newtheorem{example}[theorem]{Example}

\newcommand{\Z}{{\mathbb Z}}
\newcommand{\T}{{\mathcal T}}

\newcommand{\last}{{\texttt{last}}}
\newcommand{\supp}{{\texttt{supp}}}

\newcommand{\A}{{\mathcal{A}}}

\newcommand{\Til}{{\sf{Til}}}

\title{Enumeration of Colored Tilings on Graphs via Generating Functions}
\date{\today}

\begin{document}

\author{Jos\'e L. Ram\'{\i}rez}
\address{\noindent Departamento de Matem\'aticas, Universidad Nacional de Colombia, Bogot\'a,  COLOMBIA}
\email{jlramirezr@unal.edu.co}
\urladdr{http://sites.google.com/site/ramirezrjl}

\author{Diego Villamizar}
\address{\noindent  Department of Mathematics,   Xavier University of Louisiana,  New Orleans, LA 70125}
\email{dvillami@xula.edu}
\urladdr{https://sites.google.com/view/dvillami/}

\begin{abstract}
In this paper, we study the problem of partitioning a graph into connected and colored components called blocks. Using bivariate generating functions and combinatorial techniques, we determine the expected number of blocks when the vertices of a graph $G$, for $G$ in certain families of graphs, are colored uniformly and independently. Special emphasis is placed on graphs of the form $G \times P_n$, where $P_n$ is the path graph on $n$ vertices. This case serves as a generalization of the problem of enumerating the number of tilings of an $m \times n$ grid using colored polyominoes.
 \end{abstract}

\maketitle

\section{Introduction}

In this paper, we address the problem of enumerating the number of tilings of an $m\times n$ grid using colored polyominoes. This problem is part of a broader class of results concerning the enumerating of ways to partition a collection of objects into smaller sets according to specific rules. One of the simplest examples in this context is counting the number of tilings of a $2\times n$ rectangle using vertical and horizontal dominoes, which is given by the Fibonacci number $F_{n+1}$, see \cite{belcastro} for a generalization of this problem on surfaces. The more general problem of tiling an  $m\times n$   grid was independently solved by Temperley and Fisher  \cite{TM} and Kasteleyn \cite{Kas}. The solution is elegantly expressed by the following formula:
\begin{equation}
\prod_{j=1}^{\lceil \frac{m}{2}\rceil} \prod_{k=1}^{\lceil \frac{n}{2}\rceil}\left( 4 \cos^2\frac{j\pi}{m+1}+4\cos^2\frac{k\pi}{n+1}\right).
\end{equation}

\begin{figure}[H]
\centering
  \includegraphics[scale=1.2]{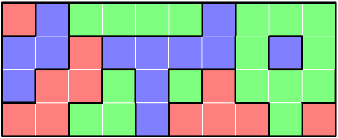}
  \caption{A tiling  of a $4 \times 10$  grid using polyominoes.}
  \label{Fig1}
\end{figure}

In the plane $\Z \times \Z$, a \emph{cell} is a unit square with vertices that have integer coordinates. A \emph{polyomino} is a finite set of cells whose interior is connected. A variation of the previous problem is to consider the number of tilings of an $m \times n$ grid using polyominoes colored with one of $k$ colors, with the condition that two adjacent polyominoes, that is, those sharing at least one edge, must have different colors. A tiling of this type is  called a \emph{$k$-colored tiling}.  
For example, Figure \ref{Fig1} shows a 3-colored tiling of a $4 \times 10$  grid. Note that this example consists of 11 polyominoes.

Richey \cite{Richey} shows that $\lim _{n,m\to \infty} e(m, n)/mn$ exists and is finite, where $e(m, n)$ is the expected number of polyominoes on an $m\times n$ grid.  Mansour \cite{Mansour}  uses automata theory to obtain a solution for grids with height at most 3.  Rolin and Ugolnikova \cite{RolUg} also apply automata theory to  the case of square polyominoes of sizes $1\times 1$ and $2 \times 2$.   Ram\'irez and Villamizar \cite{RamVil} consider this problem for square and hexagonal grids by using generating functions.  Recently, Do\v{s}lic and Podrug \cite{DP} also explored the case of hexagonal grids. Bodini \cite{Bodini} addressed a related problem, known as \emph{rectangular shape partitions}.

This combinatorial problem can be described in terms of graphs \cite{RamVil2}.  Let $G_1 = (V_1,E_1)$ and $G_2=(V_2,E_2)$ be two undirected graphs. The \emph{product} of $G_1$ and $G_2$ is defined as 
$G_1 \times G_2 = (V_1\times V_2, E_{G_1\times G_2})$, where 
\begin{multline*}
    E_{G_1\times G_2}=\{\{(v_1,v_2),(w_1,w_2)\}: (v_1=w_1 \text{ and }\{v_2,w_2\}\in E_{2})\text{ or } \\ (v_2=w_2 \text{ and }\{v_1,w_1\}\in E_{1})\}.
\end{multline*}
Let $P_n$ be a \emph{path graph}, which is a simple graph with $n$ vertices arranged in a linear sequence such that two vertices are adjacent if they are consecutive in the sequence,
and non-adjacent otherwise. A \emph{grid graph} of size $m \times n$ is defined as the product $L_{m,n}:=P_m \times P_n$.

\begin{figure}[ht]
\centering
  \includegraphics[scale=1.4]{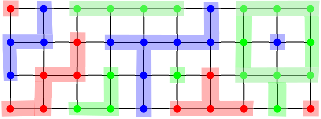}
  \caption{A 3-colored partition of size $11$ of $L_{4,10}$.}
  \label{Fig2}
\end{figure}

Let $G=(V,E)$ be an undirected graph. Two non-empty disjoint subsets $V_1, V_2\subseteq V$ are said to be \emph{neighbors} if there exists an edge $(v_1, v_2)\in E$ such that $v_1\in V_1$ and $v_2\in V_2$. A \emph{$k$-colored partition} of the vertices $V$ of $G$ is a partition  $V = \bigcup_{i=1}^s V_i$, with $s\geq 1$,  such that each $V_i$  induces a  connected subgraph of $G$, all vertices in $V_i$ are colored with  one of $k$ colors, and any pair of neighboring sets $V_i$ and $V_j$ are colored with different colors. Each set $V_i$ is called a \emph{block} of the partition, and the number of blocks is called the \emph{size} of the partition.  For example, Figure \ref{Fig2} shows a 3-colored partition of size $11$ of $L_{4,10}$. 

The problem of counting $k$-colored tilings of the grid $m \times n$  is equivalent to counting $k$-colored partitions of the graph $L_{m,n}$.  For example, the partition shown in Figure \ref{Fig2} corresponds to the 3-colored tiling given in Figure \ref{Fig1}.

Given a $k$-colored partition  $T$, we use $\rho(T)$  to denote the size of the partition $T$.  Let $G$ be an undirected graph and denote by  $\T^{(k)}(G)$  the set of all $k$-colored partitions of $G$. The random variable  $X_{\Til_{k}}(G)$  gives the size of a random $k$-colored partition in $\T^{(k)}(G)$.

\section{Colored Partitions on Trees}

The aim of this section is to enumerate the number of $k$-colored partitions of a given tree. Recall that a \emph{tree} is a connected graph with no cycles.   A \emph{perfect binary tree} is a tree in which every vertex has either $0$ or $2$ children, and all leaves are at the same height.  Let $B_n$ denote the set of perfect binary trees in which the height of each leaf is $n$. 

For fixed positive integers $n$ and $k$, we define the bivariate generating function 
 \begin{align*}
     F_{B_n}^{(k)}(x,y):=\sum_{n\geq 0}x^{n}\sum_{T\in\T^{(k)}(B_n)}  y^{\rho(T)}.
 \end{align*}

Let $F_n^{(k)}(y)$ denote the coefficient of $x^n$  in the generating function $F^{(k)}_{\mathcal{B}_n}(x,y)$, that is, $F_n^{(k)}(y)=[x^n]F^{(k)}_{\mathcal{B}_n}(x,y)$. 

\begin{lemma}
\label{lemmaArbolesBinarios}
For all $n\geq 0$, we have 
$$B_n^k(y) = ky((k-1)y+1)^{2^{n+1}-2}.$$    
\end{lemma}
\begin{proof}
We begin by noting that for $n=0$, we have $F_0^{(k)}(y) = ky$, since there is only one vertex in the tree. For $n>0$, a perfect binary tree of height $n$  consists of a root connected to two perfect binary trees of height $n-1$. Using this structure, we obtain the following recurrence relation:
    $$F_n^{(k)}(y) =F_{n-1}^{(k)}(y)^2\left ((k-1)^2y+2(k-1)+\frac{1}{y}\right), \quad n\geq 1.$$
Let $P_k(y)=(k-1)^2y+2(k-1)+1/y$. Then, we can write
    \begin{align*}
        F_n^{(k)}(y) &=F_{n-1}^{(k)}(y)^2P_k(y)=F_{n-2}^{(k)}(y)^{2^2}P_k(y)^{1+2}= \cdots = F_0^{(k)}(y)^{2^n}P_k(y)^{1+2+2^2+\cdots +2^{n-1}}.
    \end{align*}
Therefore,  we have
$$F_n^{(k)}(y)  = ky^{2^n}\left ((k-1)^2y+2(k-1)+\frac1y\right )^{2^n-1} = k\frac{y^{2^n}}{y^{2^n-1}}(((k-1)y+1)^2)^{2^n-1}.$$ 
Simplifying this expression, we obtain the desired result.
\end{proof}

\begin{figure}[H]
\centering
  \includegraphics[scale=0.8]{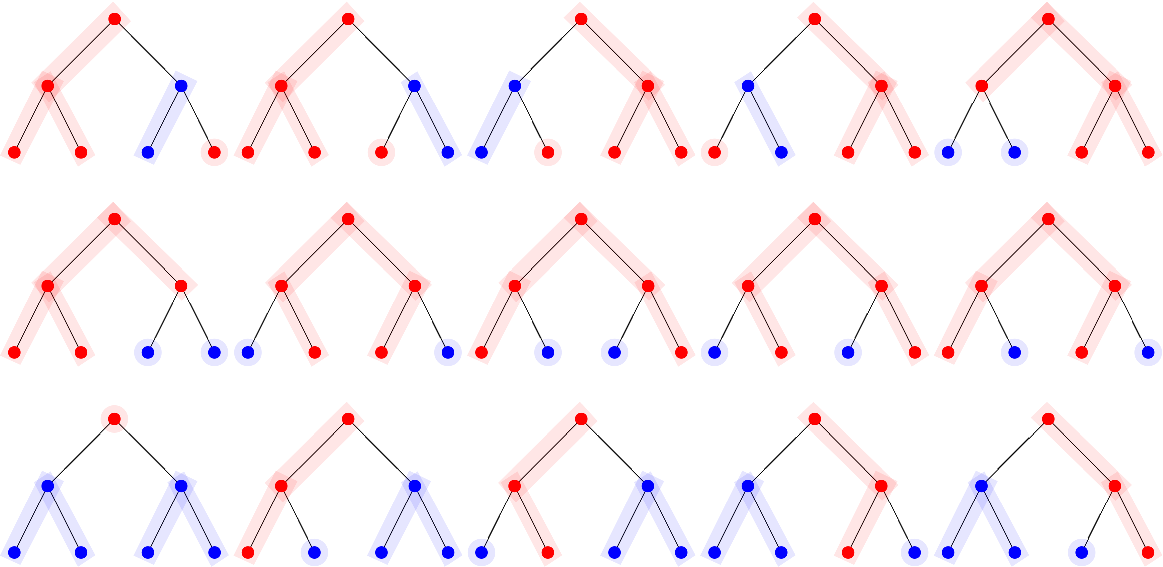}
  \caption{The 2-colored partitions of size 3 for the perfect binary tree $B_2$.}
  \label{Fig3a}
\end{figure}

For example, when $n=2$ and $k=2$, the polynomial is 
$$F_2^{(2)}(y)=2 y + 12 y^2 + \bm{30} y^3 + 40 y^4 + 30 y^5 + 12 y^6 + 2 y^7.$$
In Figure \ref{Fig3a}, we illustrate the corresponding 2-colored partitions of size 3 for the perfect binary tree $B_2$, where the root is colored red.  

From the previous lemma, we can see that
$$F_{B_n}^{(k)}(x,y)= ky+\frac{ky}{((k-1)y+1)^2}\sum _{n \geq 1}\left(((k-1)y+1)^2\right )^{2^n}x^n.$$

Using Lemma \ref{lemmaArbolesBinarios}, if we differentiate and normalize by $k^{2^{n+1}-1}$, the following corollary follows. 

\begin{corollary}
The expected number of blocks when uniformly and independently coloring a perfect binary tree with $n$ levels using $k$ colors is given by
$$\mathbb{E}[X_{\Til_k}(B_n)] = \frac{k+2(k-1)(2^{n}-1)}{k}.$$    
\end{corollary}

We can generalize this result to trees with $n$ vertices.  Let $T_n$ be a tree with $n$ vertices. We introduce the bivariate generating function
 \begin{align*}
     F_{T_n}^{(k)}(x,y):=\sum_{n\geq 0}x^{n}\sum_{T\in\T^{(k)}(T_n)}  y^{\rho(T)}.
 \end{align*}

\begin{theorem} For all $n\geq 1$, we have
$$[x^n]F_{T_n}^{(k)}(x,y) = ky((k-1)y+1)^{n-1}.$$
\end{theorem}
\begin{proof}
    To prove this result, consider an arbitrary vertex of the tree, which we designate as the  \textit{root}, and color it with one of the $k$ available colors. This contributes the term $ky$ to the generating function. Next, for each child of the root, there are two possibilities: either we color the child the same as the root, which does not create a new block, or we color the child differently. If we choose a different color, we must select one from the remaining $k-1$ colors, thus creating a new colored block. These options translate into the term $1+(k-1)y$. Repeating this process for each subsequent level down to the leaves of the tree results in the final expression.
\end{proof}

\begin{corollary}
   The expected number of blocks when uniformly and independently coloring a tree with $n$ vertices using $k$ colors is given by $$\mathbb{E}\left [X_{\Til_{k}}(T_n)\right ] = \frac{(k-1)n+1}{k}.$$
\end{corollary}
\begin{proof}
    Taking the derivative of the generating function $F_{T_n}^{(k)}(x,y)$ with respect to $y$, we obtain
    \begin{align*}
      \frac{\partial \left (ky((k-1)y+1)^{n-1}\right )}{\partial y}  &= k\left (((k-1)y+1)^{n-1}+y(n-1)(k-1)((k-1)y+1)^{n-2}\right ).
    \end{align*}
Evaluating this expression at $y=1$, we have
    $$\mathbb{E}\left [X_{\Til_{k}}(T_n)\right ] =\frac{k(k^{n-1}+(n-1)(k-1)k^{n-2})}{k^n}.$$
Simplifying the expression, we obtain the desired result.
\end{proof}

As an example, consider the case $m=1$ of the grid $L_{m,n}$, which is a tree (the path graph with $n$ vertices). In this case, we recover Theorem 1 in \cite{RamVil}.

\section{Colored Partitions on the Cycle Graph}

A \emph{walk} in a graph is a sequence of vertices $v_1, v_2, \dots,v_{\ell+1}$, not necessarily distinct, such that each pair of consecutive vertices $(v_i,v_{i+1})$ is an edge in the graph.  The \emph{length} of a walk corresponds to the number of edges, which is equivalent to the integer $\ell$.  A \emph{closed walk} is a walk in which the first vertex $v_1$ is the same as the last vertex $v_{\ell+1}$. 

Let $C_n$ denote the cycle graph on the $n$ vertices. We introduce the bivariate generating function
 \begin{align*}
     F_{C_n}^{(k)}(x,y):=\sum_{n\geq 3}x^{n}\sum_{T\in\T^{(k)}(C_n)}  y^{\rho(T)} = \sum _{n \geq 3}\sum _{i \geq 1}f_{k}(n,i)x^ny^i,
 \end{align*}
where $f_{k}(n,i)$ is the number of $k$-colored partitions of the cycle graph $C_n$ of size $i$.
\begin{theorem}
    The number of $k$-colored partitions of size $i$ for the cyclic graph $C_n$, for $n\geq 3$, is given by 
    $$f_{k}(n,i)=\begin{cases}
        k, & \text{ if }i=1;\\
        \binom{n}{i}\left ((k-1)^i+(k-1)(-1)^i\right), & \text{otherwise.}
    \end{cases}$$
 Moreover, its bivariate generating function $F_{C_n}^{(k)}(x,y)$ is given by 
    $$\frac{ky}{1-x}+\frac{1}{1-x(1+(k-1)y)}+\frac{k-1}{1-x(1-y)}-\frac{k}{1-x}-(1+x)ky-(ky+k(k-1)y^2)x^2.$$
\end{theorem}
\begin{proof}
    If the size of the partition is one, then all vertices are colored the same color. Consequently,  we can choose this color in $k$ ways. For $i>1$, consider the choice of colors as walks in a complete graph with $k$ vertices, where each vertex is labeled with one of the $k$ colors. The walk starts and ends at either the same color (which means that the size of the partition decreases by one) or at different colors. Recall the known formula (cf.\  \cite[p. 5]{Stanley}) for the number of closed walks of size $\ell$ in the complete graph $K_m$, given by 
$$C_{K_m,\ell}=\frac{1}{m}(\left (m-1)^{\ell} +(m-1)(-1)^{\ell}\right ),$$
and the number of walks that start and end at different vertices is given by 
$$W_{K_m,\ell}=\frac{1}{m}\left ((m-1)^{\ell}-(-1)^{\ell}\right ).$$
The number of ways to select the blocks is given by the number of integer compositions of $n$ into $i$ parts, which is $\binom{n-1}{i-1}$. Depending on the colors we choose, there will either be $i$ blocks or $i-1$ blocks.  If the walk starts and ends in different colors, then we have $i$ blocks. If the walk starts and ends at the same color, then we obtain $i-1$ blocks. Thus, by considering both cases for the starting and ending colors,  the number of $k$-colored partitions of size $i$ is given by
\begin{align*}
    f_{k}(n,i)&= 2\binom{n-1}{i-1}\binom{k}{2}W_{K_k,i-1}+\binom{n-1}{i}\binom{k}{1}C_{K_k,i}\\
    &=2\binom{n-1}{i-1}\binom{k}{2}\frac{1}{k}\left ((k-1)^{i-1}-(-1)^{i-1}\right )+\binom{n-1}{i}\binom{k}{1}\frac{1}{k}\left ((k-1)^i+(k-1)(-1)^i\right )\\
    &=\binom{n-1}{i-1}\left ((k-1)^{i}+(k-1)(-1)^{i}\right )+\binom{n-1}{i}\left ((k-1)^i+(k-1)(-1)^i\right )\\
    &=\binom{n}{i}\left ((k-1)^i+(k-1)(-1)^i\right ).
\end{align*}
The last equality is Pascal's recursion for the binomial coefficients. Finally, using the expression for the number of $k$-colored partitions, we obtain the bivariate generating function
\begin{align}
\label{cycleGF}
    F_{C_n}^{(k)}(x,y)&= \sum _{n\geq 3}x^n\left (ky+\sum _{i\geq 2}\binom{n}{i}\left ((k-1)^i+(k-1)(-1)^i\right )y^i\right )\\
    &=\sum _{n\geq 3}x^n\left (ky+((k-1)y+1)^n+(k-1)(1-y)^n-k\right ).
\end{align}
The result follows by completing the geometric series.
\end{proof}

\begin{figure}[H]
\centering
  \includegraphics[scale=0.8]{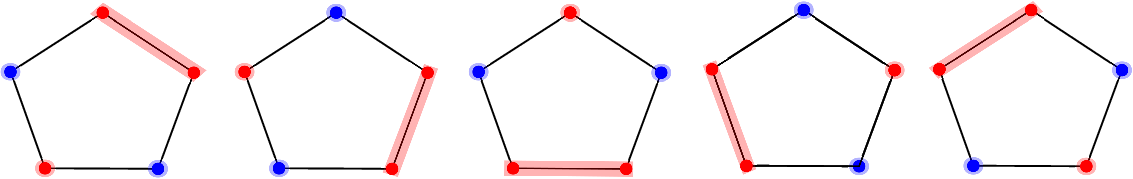}
  \caption{The 2-colored partitions of size 4 for the cycle graph $C_5$.}
  \label{Fig4}
\end{figure}

For example, when $n=5, i=4$, and $k=2$, we have $f_2(5,2)=10$.  In Figure \ref{Fig4}, we illustrate the corresponding 2-colored partitions of size 4 for the cycle graph $C_5$, where two consecutive vertices are colored red.  

\begin{corollary}
   The expected number of blocks when uniformly and independently coloring a cycle graph with $n$ vertices using $k$ colors is given by
   \begin{align*}\mathbb{E}[X_{\Til_k}(C_n)] = \frac{1}{k^n}\left (k+n(k^n-k^{n-1})\right ).\end{align*}
\end{corollary}
\begin{proof}
 Starting with the expression of the generating function in \eqref{cycleGF}, we differentiate with respect to $y$ to obtain 
    $$k+n(k-1)((k-1)y+1)^{n-1}-n(k-1)(1-y)^{n-1}.$$ Evaluating this derivative at $y=1$ gives the desired result.
    \end{proof}

\section{Colored Partitions on the Complete Graph}

Let $K_n$ denote the complete graph on $n$ vertices. We introduce the exponential bivariate generating function
 \begin{align*}
     G_{K_n}^{(k)}(x,y):=\sum_{n\geq 1}\frac{x^{n}}{n!}\sum_{T\in\T^{(k)}(K_n)}  y^{\rho(T)} = \sum _{n \geq 1}\sum _{i\geq 1}g_{k}(n,i)\frac{x^n}{n!}y^i,
 \end{align*}
 where $g_{k}(n,i)$ denotes the number of $k$-colored partitions of size $i$ on the complete graph $K_n$.

\begin{theorem}
\label{thmComplete}
    The number of $k$-colored partitions of size $i$ for the complete graph $K_n$, for $n\geq 1$, is given by 
    $$g_{k}(n,i)={n\brace i}\binom{k}{i}i!.$$
Moreover, its exponential bivariate generating function is given by 
    $$G_{K_n}^{(k)}(x,y)=(1+y(e^x-1))^k.$$
\end{theorem}
\begin{proof}
    To construct a $k$-colored partition of size $i$, we first choose $i$ colors from the available $k$ colors, which can be done in $\binom{k}{i}$ ways. Next,  we partition the $n$ vertices of the complete graph $K_n$ into $i$ disjoint non-empty blocks and assign one of chosen  colors to each block.  This can be done in ${n\brace i}i!$ ways, where ${n \brace i}$ are the Stirling numbers of the second kind. The assignment of the colors gives the $k$-colored partitions because the graph is complete.
    
    Using the above, we have that the generating function is given by
    \begin{align*}
        G_{K_n}^{(k)}(x,y) &=\sum _{n\geq 0}\frac{x^n}{n!}\sum _{i = 0}^n {n\brace i}\binom{k}{i}i!y^i\\
        &=\sum _{i\geq 0}\binom{k}{i}i!y^i\sum _{n\geq i}{n\brace i}\frac{x^n}{n!}\\
        &=\sum _{i\geq 0}\binom{k}{i}i!y^i\frac{(e^x-1)^i}{i!}\\
        &=\left (y(e^x-1)+1\right )^k. & \qedhere
    \end{align*}
\end{proof}

\begin{figure}[H]
\centering
  \includegraphics[scale=0.85]{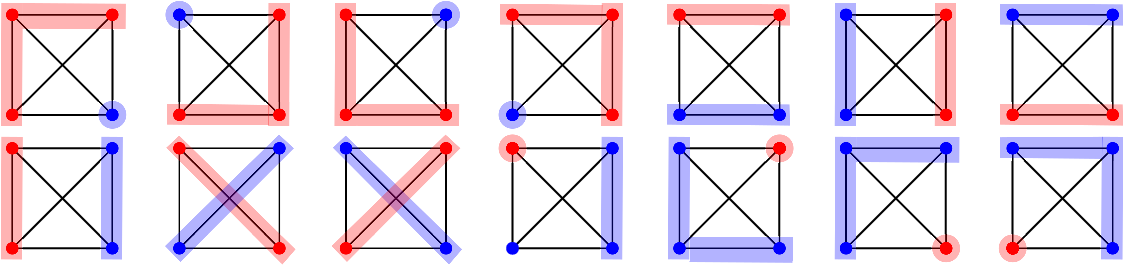}
  \caption{The 2-colored partitions of size 2 for the complete graph $K_4$.}
  \label{Fig5}
\end{figure}

For example, when $n=5, i=4$, and $k=2$, we have $g_2(4,2)=14$.  In Figure \ref{Fig5}, we illustrate the corresponding 2-colored partitions of size 2 for the complete graph $K_4$.  

\begin{corollary}
\label{crlComplete}
   The expected number of blocks when uniformly and independently coloring a complete graph with $n$ vertices using $k$ colors is given by $$\mathbb{E}\left [X_{\Til_k}(K_n)\right ]=k-\frac{(k-1)^n}{k^{n-1}}.$$
\end{corollary}
\begin{proof}
    Using Theorem \ref{thmComplete}, to compute the expected value, we need to find the coefficient of $x^n/n!$ in the expansion of  $$ \frac{\partial G_{K_n}^{(k)}(x,y)}{\partial y} = k(e^x-1)\left (y(e^x-1)+1\right )^{k-1}.$$ Evaluating at $y=1$, we have $(k(e^x-1)(e^{x(k-1)})$. Extracting the coefficient of $x^n$ from the expansion, we obtain
    $$[x^n](k(e^{kx}-e^{(k-1)x}) = \frac{k(k^n-(k-1)^n)}{n!}.$$
    Dividing by the total number of colorings, $k^n$, we have the expected value \[\mathbb{E}\left [X_{\Til_k}(K_n)\right ] = \frac{k(k^n-(k-1)^n)}{k^n} = k-\frac{(k-1)^n}{k^{n-1}}.\qedhere  \]
\end{proof}

\section{Colored Partitions on the Complete Bipartite Graph}\label{sectionBipartite}

Let $K_{n,m}$ denote the complete bipartite graph on $n$ and $m$ vertices.

\begin{theorem}
   The expected number of blocks when uniformly and independently coloring a complete bipartite graph $K_{n,m}$ using $k$ colors is given by
    $$\mathbb{E}[X_{\Til_k}(K_{n,m})] =\frac{nk^n(k-1)^m+mk^m(k-1)^n+k(k^n-(k-1)^n)(k^m-(k-1)^m)}{k^{n+m}}.$$
\end{theorem}
\begin{proof}
A colored partition is formed by assigning a color to each vertex of a bipartite graph. If the same color $x\in [k]$ is used on the vertices of both parts, they will merge into a single block. Therefore, the number of blocks is equal to the total number of vertices $n+m$ minus the number of vertices that share the same color in both parts. This implies that the expected number of blocks is
\begin{align}
\label{eqBip}
    \mathbb{E}[X_{\Til_k}(K_{n,m})] = \frac{1}{k^{n+m}}\sum _{\substack{f :[n]\rightarrow [k]\\g:[m]\rightarrow [k]}}\left ( n+m-\sum _{x\in \texttt{Im}(f)\cap \texttt{Im}(g)}\left (|f^{-1}(x)|+|g^{-1}(x)|-1\right )\right ).
\end{align}
    We will divide the sum into two parts:
    \begin{enumerate}
     \item First, consider $$\sum_{\substack{f :[n]\rightarrow [k]\\g:[m]\rightarrow [k]}}\sum _{\texttt{Im}(f)\cap \texttt{Im}(g)}|f^{-1}(x)|.$$
        This summation is equivalent to iterating over each element on $[n]$, one block of the graph, and counting the pairs of functions where $f(x)$ is in the image of $g$. By symmetry, we find that this summation is equal to 
        $$n|\{(f,g): f(1)\in \texttt{Im}(g)\}|.$$ To compute the size of this set, we iterate over the size of the image of $g$ and use the fact that a function is surjective over its image. Given the size of the image, we choose an element of the image,  leading to the result
        $$|\{(f,g): f(1)\in \texttt{Im}(g)\}| = \sum _{\ell = 1}^k \ell \binom{k}{\ell}{m \brace \ell}\ell ! k^{n-1}.$$
        Using the computation from the proof of Corollary \ref{crlComplete}, we obtain that this summation is equal to $$nk^{n-1}\left ((k+1)k^m-k(k-1)^m-k^m\right ) = n\left (k^{m+n}-k^n(k-1)^m\right).$$
        
    \item Next, consider $$\sum _{\substack{f :[n]\rightarrow [k]\\g:[m]\rightarrow [k]}}|\texttt{Im}(f)\cap \texttt{Im}(g)|.$$
        This summation is equivalent to iterating over the $k$ colors and counting the number of function pairs  that include  that color in their image. By symmetry, we find that the summation is equal to
        $$k|\{(f,g): 1\in \texttt{Im}(f)\cap \texttt{Im}(g)\}| = k(k^n-(k-1)^n)(k^m-(k-1)^m).$$
            \end{enumerate}
 Substituting these results into the expression in \eqref{eqBip} gives the result.   
\end{proof}

\section{Colored Partitions in the Cartesian Product}
\begin{figure}[ht]
\centering
  \includegraphics[scale=1.2]{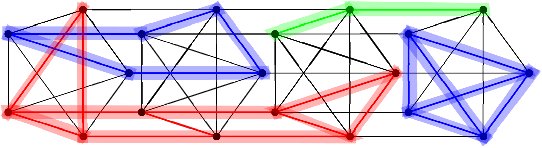}
  \caption{A 3-colored partition of size $4$ of $K_5\times P_4$.}
  \label{Fig3}
\end{figure}

In this section, we enumerate the $k$-colored partitions for the graph product $K_m \times P_n$ for $m=3, 4$.  For example,  Figure \ref{Fig3} shows a 3-colored partition of size $4$ for the graph $K_5\times P_4$.

\subsection{The case $m=3$.}
In this section, we provide the explicit bivariate generating function for the $k$-colored partitions of  $U_n:=K_3 \times P_n$ for all $n\geq 1$.
For fixed positive integers $n$ and $k$, we define the bivariate generating function as 
 \begin{align*}
    T(x,y):=\sum_{n\geq 0}x^{n}\sum_{T\in\T^{(k)}(U_n)}  y^{\rho(T)}, 
 \end{align*}
 where $\T^{(k)}(U_n)$ denotes the set of $k$-colored partitions of $U_n$, and $\rho(T)$ is the size of the partition $T$.
\begin{theorem}\label{teo2col}
The bivariate generating function $T(x,y)$ is given by $p(x,y)/q(x,y)$, where 
\begin{multline*}
   p(x,y)= k x y \left(1 - (1 - k) y (3 - (2 - k) y) - 
    x (1 - y) (4 - y (13 - 5 k) + (3 + k) (3 - 2 k) y^2 \right. \\
\left.    - (1 - k) (1 + (3 - k) k) y^3) +  x^2 (1 - y)^2 (3 - (2 + 4 k) y + (1 - k + 3 k^2) y^2 + 
       k^2 (1 - k ) y^3)\right)
\end{multline*}
and \footnotesize
\begin{multline*}
   q(x,y)= 1 - x \left(5 - 12 (2 - k) y + 
     x^2 (1 - y)^2 (3 - (2 + 4 k) y + k (1 + 2 k) y^2 + 
        k (1 - k ) y^3 - (k - 1)^4 y^4)   \right. \\ \left.
-     x (1 - y) (7 - (25 - 8 k) y + (20 - 3 k - 4 k^2) y^2 + (5 -   24 k + 21 k^2 - 5 k^3) y^3 - (7 - 18 k + 17 k^2 - 7 k^3 + k^4) y^4)  \right. \\ \left.
    + y^2 (32 - 26 k + 6 k^2 - (13 - 14 k + 6 k^2 - k^3) y)\right).
\end{multline*}
\normalsize 
Moreover, $[x^n]T(x,1)=k^{3n}$.
\end{theorem}
\begin{proof}
Let $\mathcal{A}_{n, i}$ denote the set of $k$-colored tilings in $\T^{(k)}(U_n)$ such that the last triangle is colored  with exactly $i$ different colors, where $i=1, 2, 3$. For a $k$-colored partition  $T\in\T^{(k)}(U_n)$, we denote the last triangle of $T$ by $\last(T)$.

Define the bivariate generating functions
 \[T_i(x,y):=\sum_{n\geq 1}x^{n}\sum_{T\in\mathcal{A}_{n,i}}  y^{\rho(T)}, \quad \text{for } i=1, 2, 3.\] 
It is clear that $T(x,y)=T_1(x,y)+T_2(x,y) + T_3(x,y)$.  

Let $T$ be a $k$-colored partition in $\mathcal{A}_{n,1}$. If $n=1$, then $T=K_3$, and its contribution  to the generating function is the term $kxy$, as it must be monochromatic. If $n>1$, then $T$  can be decomposed as $T_i K_3$, where $K_3=\last(T)$ is monochromatic, and $T_i\in \mathcal{A}_{n-1,i}$, for $i=1, 2, 3$.  
The decomposition leads to:
\begin{itemize}
    \item Case 1 ($i=1$).  In this case, $\last(T_1)$ and $K_3$ may or may not share the same color. Therefore, the generating function is given by
     $$xT_1(x,y) + (k-1)xyT_1(x,y).$$
    \item Case 2 ($i=2$). Here there are two possibilities, $\last(T_2)$ and $K_3$ either coincide in one of the colors of $\last(T_2)$ or are colored differently. Thus, the generating function is
       $$2xT_2(x,y) + (k-2)xyT_2(x,y).$$
       
    \item Case 3 ($i=3$).  Similarly, $\last(T_3)$ and $K_3$ may coincide in one of the colors of $\last(T_3)$ or  be entirely different. The generating function for this case is
      $$3xT_3(x,y) +(k-3)xyT_3(x,y).$$
\end{itemize}
Table \ref{deco1} illustrates the three cases discussed above.
 \begin{table}[ht]
        \centering 
        \begin{tabular}{|c|}  \hline
        Case 1\\
\includegraphics[scale=0.9]{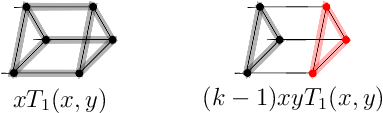}   \\ \hline \hline   Case 2\\
\includegraphics[scale=0.9]{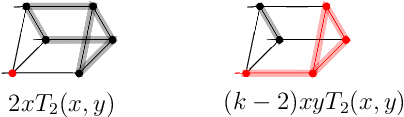}  \\ \hline \hline   Case 3\\
\includegraphics[scale=0.9] {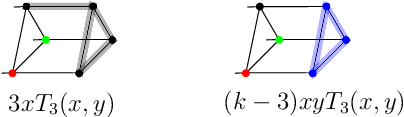}  \\ \hline    
        \end{tabular}
        \caption{Cases for the generating function $T_1(x,y)$.}
        \label{deco1}
    \end{table}
    
    Combining these,  we obtain the following functional equation:    
\begin{multline*}
    T_1(x,y) = kxy + xT_1(x,y) + (k-1)xyT_1(x,y) +2xT_2(x,y)\\ + (k-2)xyT_2(x,y) + 3xT_3(x,y) +(k-3)xyT_3(x,y).
\end{multline*} 

For $T\in\mathcal{A}_{n,2}$, we  analyze its decomposition as shown in  Table \ref{deco2}. Each case yields terms based on the interaction of $\last(T)$ with $K_3$. 

    \begin{table}[H]
        \centering 
        \begin{tabular}{|c|}
        \hline Case 1\\
  \includegraphics[scale=0.7] {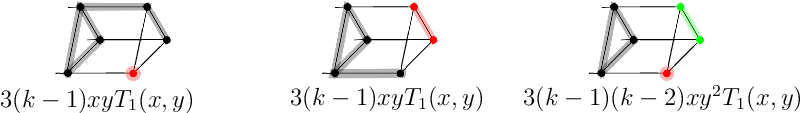}    \\ \hline \hline
               \hline Case 2\\
  \includegraphics[scale=0.7] {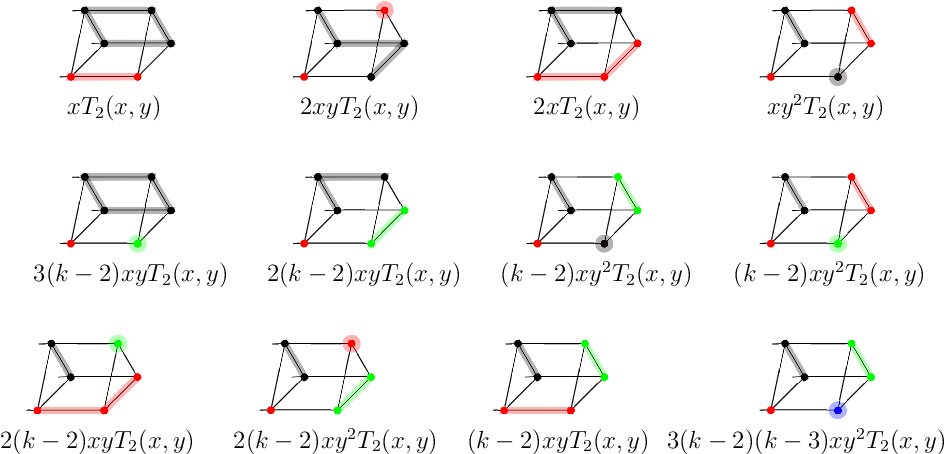}    \\ \hline 
   \hline Case 3\\
  \includegraphics[scale=0.7] {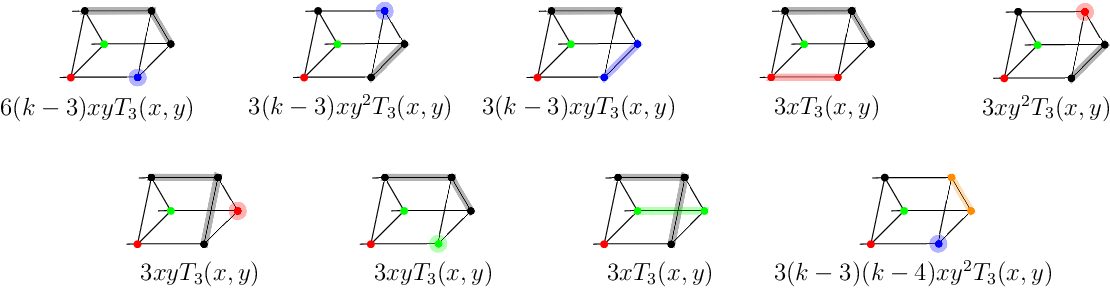}    \\ \hline 
          \end{tabular}
        \caption{Cases for the generating function $T_2(x,y).$
        }
        \label{deco2}
    \end{table}
   
From this decomposition the generating function for $T_2(x,y)$ satisfies:
\begin{multline*}
    T_2(x,y)= 3k(k - 1)xy^2 + 6(k - 1)xyT_1(x,y) + 
 3(k-1)(k - 2)xy^2T_1(x,y) \\
 + (3x+2xy+xy^2)T_2(x,y) + (8 (k - 2)xy + 4(k-2)xy^2) T_2(x,y) \\ + 
 3 (k - 2)(k - 3) xy^2T_2(x,y) + (9(k - 3)xy+ 9(k - 3) xy^2)T_3(x,y) \\ + (6 x + 6 x y + 6 x y^2)T_3(x,y)  + 3 (k - 3) (k - 2) x y^2T_3(x,y).
\end{multline*}

Finally, let $T$ be a $k$-colored partition in $\mathcal{A}_{n,3}$. Again, from a similar argument as before, we can obtain the different  possible decompositions, as show in Table \ref{deco3}.

    \begin{table}[H]
        \centering 
        \begin{tabular}{|c|}
        \hline Case 1\\
  \includegraphics[scale=0.7] {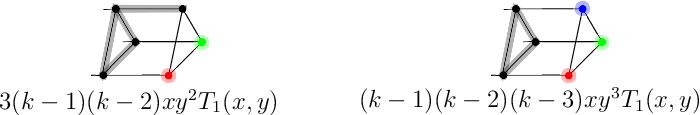}    \\ \hline \hline
               \hline Case 2\\
  \includegraphics[scale=0.7] {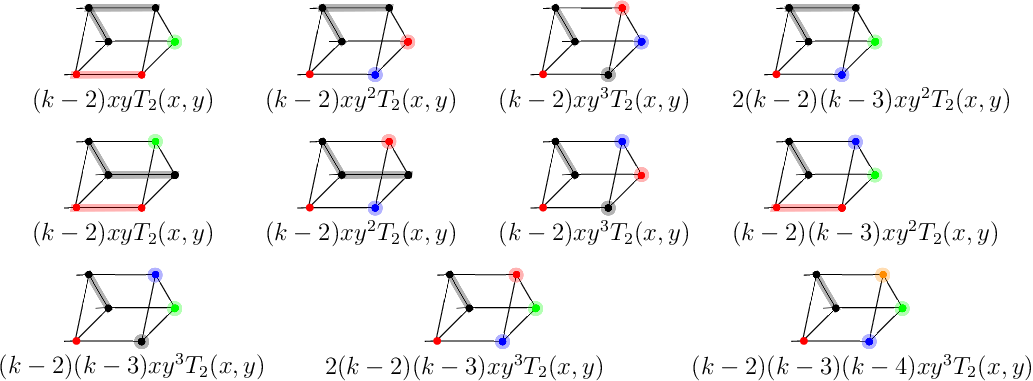}    \\ \hline 
   \hline Case 3\\
  \includegraphics[scale=0.7] {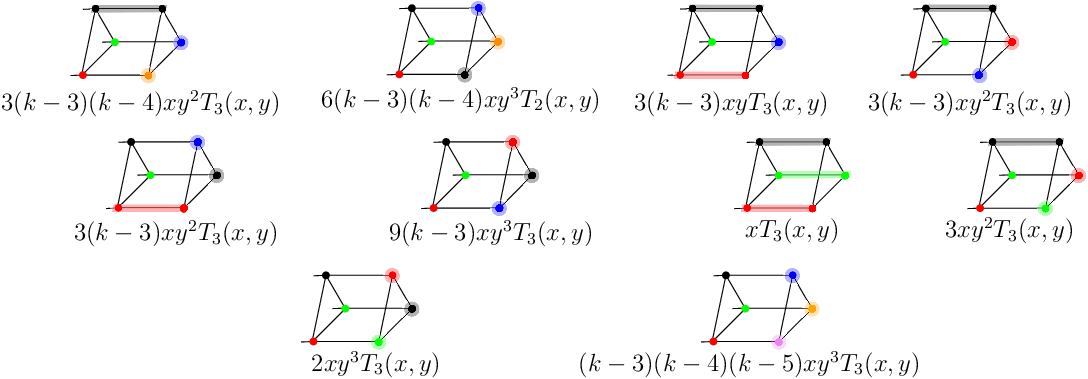}    \\ \hline 
          \end{tabular}
        \caption{Cases for the generating function $T_3(x,y).$
        }
        \label{deco3}
    \end{table}
From this decomposition we obtain the functional equation:
\begin{multline*}
    T_3(x,y)= k(k-1)(k-2)xy^3 + 
 3(k - 1)(k - 2)xy^2 T_1(x,y) + (k - 1) (k - 2)(k - 3) x y^3 T_1(x,y) \\+ (2 (k - 2) x y + 2 (k - 2) x y^2 + 
    2 (k - 2) x y^3) T_2(x,y) \\ + (3 (k - 2) (k - 3) x y^2 + 
    3 (k - 2) (k - 3) x y^3) T_2(x,y)  + (k - 2) (k - 3) (k - 4) xy^3T_2(x,y)\\ + 
 3 (k - 3) (k - 4) x y^2T_3(x,y) + 6 (k - 3) (k - 4) x y^3 T_3(x,y) + 
 3 (k - 3) x y T_3(x,y) \\ + 6 (k - 3) x y^2 T_3(x,y) + 
 9 (k - 3) x y^3 T_3(x,y)  + (x + 3 x y^2 + 2 x y^3) T_3(x,y) \\ + (k - 3) (k - 4) (k - 5) xy^3 T_3(x,y).
\end{multline*}
    Since $T(x,y)=T_1(x,y)+T_2(x,y)+T_3(x,y)$, we obtain a system of four linear equations with four unknowns $T(x,y), T_1(x,y), T_2(x,y)$, and $T_3(x,y)$. Solving the system for $T(x,y)$ yields the desired result. 
\end{proof}

\begin{corollary}
   The expected number of blocks when uniformly and independently coloring the graph $U_n=K_3 \times P_n$ with $3n$ vertices using $k$ colors is given by
    $$\mathbb{E}\left [X_{\Til_k}(U_n)\right ]=\frac{k^{3n-5}((1 - 3 k^2 + 3 k^4)  -(1 - k)^3 (1 + 3 k + 3 k^2) n)}{k^{3n}}.$$
\end{corollary}

For example, for $k=2$ we obtain the bivariate generating function (see Theorem 3.1 of \cite{RamVil2}): 
\begin{align*}
T(x,y)=\frac{2 x y (1 + 3 y - x (3 - 7 y + 4 y^2))}{1 - x (4 + 3 y + y^2) + x^2 (3 - 7 y + 3 y^2 + y^3)}.
\end{align*}
Moreover, 
    $$\mathbb{E}\left [X_{\Til_2}(U_n)\right ]=\frac{2^{3n-5}(37+19n)}{2^{3n}}.$$

\subsection{The general case}\label{gcase}

In this section, we present a general approach for counting  the number of $k$-colored partitions of the graph $U_n^{(m)}:=K_m\times P_n$. Let $\mathcal{T}^{(k)}(U_n^{(m)} )$ denote the set of possible $k$-colored partitions for $U_n^{(m)}$  For fixed positive  integers $m$ and $k$, we introduce the following bivariate generating function:  
\begin{align}
  T_{m}^{(k)}(x,y)=\sum _{n\geq 1}x^n\sum _{T\in \mathcal{T}^{(k)}(U_n^{(m)} )}y^{\rho (T)}.
\end{align}
To determine the generating function $T_{m}^{(k)}(x,y)$, we establish a system of equations indexed by the possible $k$-colored partitions of the complete graph $K_m$.

\begin{definition}
    Let $\mathcal{C}_{m,k}$ denote the set of sequences of subsets of $[m]$ that are pairwise disjoint and whose union is $[m]$. Formally,  
    $$\mathcal{C}_{m,k} =\left\{\mathcal{A}= (A_1,\dots,A_k): \bigcup_{i=1}^k A_i = [m] \text{ and }A_i\cap A_i = \emptyset\right\}.$$
\end{definition}

First, notice that these sequences can correspond to any configuration in 
 $\mathcal{T}^{(k)}(U_1^{(m)} )$ and they are in bijection with functions from $[m]$ to $[k]$ by considering the preimages. For a given $\mathcal{A}\in \mathcal{C}_{m,k}$, we denote by $\supp(\mathcal{A})$ the number of indices $i\in [k]$ for which $A_i\neq \emptyset$, these are the colors that we have used at least once on a vertex of the corresponding configuration. For  fixed integers $m$ and $k$, and a given  $\mathcal{A}\in \mathcal{C}_{m,k}$,  consider the following bivariate function:
$$T_{\mathcal{A}}(x,y) = \sum _{n\geq 1}x^n\sum _{T\in \mathcal{T}_{\mathcal{A}}^{(k)}(U_n^{(m)} )}y^{\rho (T)},$$
where $\mathcal{T}_{\mathcal{A}}^{(k)}(U_n^{(m)} )$ are the $k$-colored partitions of  $U_n^{(m)}$ that end in the configuration specified by $\mathcal{A}$.  In this way, we would have
$$T_m^{(k)}(x,y) = \sum _{\mathcal{A}\in \mathcal{C}_{m,k}}T_{\mathcal{A}}(x,y).$$

We have now the following system of equations:
$$T_{\A}(x,y) = xy^{|\supp(\A)|} + x\sum _{\mathcal{B} \in \mathcal{C}_{m,k}}y^{\left |\{i\in [k]:A_i \neq \emptyset \text{ and } A_i\cap B_i=\emptyset\}\right |}T_{\mathcal{B}}(x,y).$$

Consider, further, the following equivalence relation over the set $\mathcal{C}_{m,k}$, defined by $\mathcal{A}\sim \mathcal{B}$ if and only if there exists a permutation $\sigma \in \mathfrak{S}_k$, where $\mathfrak{S}_k$ denotes the group of permutations on $k$ elements, such that $|A_i|=|B_{\sigma (i)}|$ for all $1\leq i\leq k$.

Notice that $|\mathcal{C}_{m,k}|=k^m$, since these are all possible colorings for the complete graph. However, the number of distinct equivalence classes in $\mathcal{C}_{m,k}/\sim$ is given by the coefficient of $q^m$ in the Gaussian binomial coefficient $\binom{m+k}{k}_q$.  This is because the coefficient of $q^m$ in $\binom{m+k}{k}_q$ counts the number of partitions of $m$ into $k$ or fewer parts, where each part is bounded above by $m$. 
    
The size of each equivalence class $[\mathcal{A}]\in \mathcal{C}_{m,k}/\sim$ is given by  $$|[\mathcal{A}]|=\frac{m!}{\prod _{i\in [k]}|A_i|!}\cdot \frac{k!}{\prod _{i=0}^m\texttt{m}(\mathcal{A})_i!}.$$ 
The first multinomial considers the choice of each set based on cardinality, while the second multinomial considers a permutation on the sequence of cardinals. Here, $\texttt{m}(\mathcal{A})_i$ denotes the number of sets in the sequence $\mathcal{A}$ that have size $i$.

\begin{example}\label{ExA}
    Consider $m = 4$ and $k = 2$, then there are $16$ possible colorings listed below as elements in $\mathcal{C}_{4,2}$
    \begin{align*}
        \mathcal{C}_{4,2} =  \{&(\emptyset ,\{1,2,3,4\}),(\{1,2,3,4\},\emptyset),(\{1\},\{2,3,4\}),(\{2\},\{1,3,4\}),\\
        &(\{3\},\{1,2,4\}),(\{4\},\{1,2,3\}),(\{1,2,3\},\{4\}),(\{1,2,4\},\{3\}),\\
        &(\{1,3,4\},\{2\}),(\{2,3,4\},\{1\}),(\{1,2\},\{3,4\}),(\{1,3\},\{2,4\})\\
        &(\{1,4\},\{2,3\}),(\{2,3\},\{1,4\}), (\{2,4\},\{1,3\}),(\{3,4\},\{1,2\})\}.
    \end{align*}
The quotient is given by $\mathcal{C}_{4,2}/\sim = \{[(0,4)],[(1,3)],[(2,2)]\}$, where the tuple $(a,b)$ represents the element $(\{1,2,\dots ,a\},\{a+1,\dots, 4\})$ in $\mathcal{C}_{4,2}$. The sizes of the classes are given by 
    \begin{align*}
      |[(0,4)]| &= \frac{2!}{1!\cdot 0!\cdot 0!\cdot 0!\cdot 1!}\frac{4!}{4!\cdot 0!} =2,\\
      |[(1,3)]| &= \frac{2!}{0!\cdot 1!\cdot 0!\cdot 1!\cdot 0!}\frac{4!}{3!\cdot 1!} =8,\\
      |[(2,2)]| &= \frac{2!}{0!\cdot 0!\cdot 2!\cdot 0!\cdot 0!}\frac{4!}{2!\cdot 2!} =6.
    \end{align*}
\end{example}

This equivalence relation enables us to reduce the number of variables in the system of equations used to compute $T_m^{(k)}(x,y)$.  This reduction is possible due to the following proposition, which stems from the fact that the equivalence relation captures both the symmetries of the complete graph and the symmetries of the color permutations.
\begin{proposition}
    If $\mathcal{A}\sim \mathcal{B}$, then $T_{\mathcal{A}}(x,y)=T_{\mathcal{B}}(x,y)$.
\end{proposition}

Using the proposition above, and for $k = 2$, we can express the generating functions in terms of the tuples $(a,b)$, such that $T_{(a,b)}(x,y) = T_{\mathcal{A}}(x,y)$ for any $\mathcal{A}\in [(a,b)]$. As an example of this reduction, we compute the generating function for the number of $2$-colored partitions of the graph $K_4\times P_n$.

\begin{theorem}
    The bivariate generating function $T_4^{(2)}(x,y)$ is given by 
    \begin{align*}
      \frac{2xy\left (1+7y-x(y-1)\left (7y^2+y-9\right )+x^2(y-1)^2\left (8y^2-17y+8\right )\right )}{1-2x(y^2+2y+5)+x^2(y-1)(y^3+6y^2+8y-17)-x^3(y-1)^2(y^3+6y^2-17y+8)}.  
    \end{align*}
\end{theorem}
\begin{proof}
    Using the description above and Example \ref{ExA}, we consider three distinct variables $T_{(0,4)}(x,y)$, $T_{(1,3)}(x,y)$, and $T_{(2,2)}(x,y)$.  Each variable corresponds to a different configuration for the last copy of $K_4$, and each contributes differently to the bivariate generating function. The contribution of each variable is determined by the size of the corresponding equivalence class.   Therefore, we have
    \begin{align}
    \label{expK4}
        T^{(2)}_{4}(x,y)=2T_{(0,4)}(x,y) + 8T_{(1,3)}(x,y) + 6T_{(2,2)}(x,y).
    \end{align}
The following relations arise when considering each possible pair of contiguous copies of $K_4$:
    \begin{align*}
        T_{(0,4)}(x,y) &= xy+x(1+y)T_{(0,4)}(x,y)+8xT_{(1,3)}(x,y)+6xT_{(2,2)}(x,y),\\
      T_{(1,3)}(x,y) &= xy^2+2xyT_{(0,4)}(x,y)+x(y^2+3y+4)T_{(1,3)}(x,y)+3x(1+y)T_{(2,2)}(x,y),\\
      T_{(2,2)}(x,y) &= xy^2+2xyT_{(0,4)}(x,y)+4x(y+1)T_{(1,3)}(x,y)+x(y^2+5)T_{(2,2)}(x,y).
    \end{align*}
    As an example, consider $T_{(1,3)}(x,y)$. Figure \ref{fig:K4} illustrates all eight possible combinations that can result in a configuration ending with $(1,3)$ after a previous configuration that also ends in $(1,3)$. By solving this $3\times 3$ system of equations and substituting the solution into \eqref{expK4}, the theorem follows.
    \begin{figure}[ht!]
        \centering
        \includegraphics[scale = 0.8]{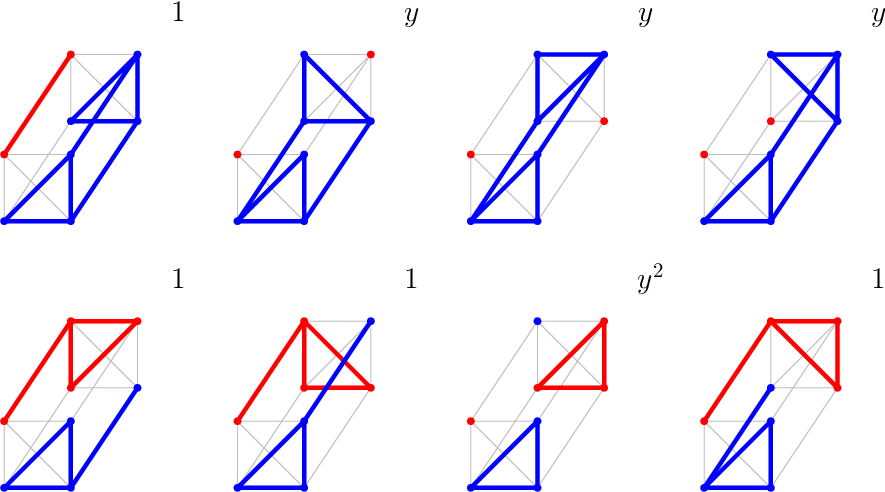}
        \caption{The contribution $y^2+3y+4$ on $T_{(1,3)}(x,y)$.}
        \label{fig:K4}
    \end{figure}

\end{proof}

\begin{corollary}\label{expval1}
   The expected number of blocks when uniformly and independently coloring the graph $K_{4}\times P_n$ with $4n$ vertices using two colors is given by
    $$\mathbb{E}\left [X_{\Til_2}(K_{4}\times P_n)\right ]=\frac{2^{4 n-7}(175 + 65 n)}{2^{4n}}.$$
\end{corollary}

Using this approach, we wrote a Python code\footnote{
\url{https://github.com/Phicar/ColoredColumns/blob/main/CompleteGraphGenFunc.py}} to calculate the generating functions. For example, this method produces the following generating functions:
\begin{align}
\resizebox{\textwidth}{!}{
\begin{minipage}{\textwidth}
    \begin{align*}
       T_5^{(2)}(x,y)=&\frac{2xy\left (1+15y-x(y-1)\left (15y^2-13y-21\right )+x^2(y-1)^2\left (16y^2-51y+30\right )\right )}{1-2x(y^2+4y+11)+x^2(y-1)(y^3+10y^2+2y-51)-x^3(y-1)^2(y^3+10y^2-51y+30)},\\
T_6^{(2)}(x,y)=&\frac{2xy\left (1+31y-x(y-1)(62y^2-103y-48)+x^2(y-1)^2(31y^3-72y^2-125y+155)-x^3(y-1)^3(32y^3-185y^2+263y-108)\right)}{1-x(3y^2+12y+49)+x^2(y-1)(3y^3+28y^2-6y-203)-x^3(y-1)^2(y^4+16y^3-40y^2-262y+263)+x^4(y-1)^3(y^4+15y^3-167y^2+263y-108)},\\
T^{(3)}_{4}(x,y)=&\frac{3xy(1+6y+2y^2+2x(y-1)(y^3-9y^2+y+2)-x^2(2y-1)(y-1)^2(9y^2-8y+3))}{1-x(2y^3+8y^2+12y+5)-x^2(y-1)(2y^4-13y^3-25y^2-y+7)+x^3(y-1)^2(16y^4+6y^3-21y^2+14y-3)}.
    \end{align*}
\end{minipage}
}
\end{align}

The  results above  lead to the following expression for the expected number of blocks in a random $k$-colored partition of the graph  $K_{m}\times P_n$.

\begin{theorem}
The expected size of a $k$-colored partition on the graph $K_{\ell}\times P_n$ assuming each vertex is  colored uniformly and independently,  is given by
\begin{align*}
    \mathbb{E}[X_{\Til_k}(K_{\ell}\times P_n)]
    &=\frac{k^{\ell n-(2\ell -1)}\left (\left (k^{2\ell}-(k^2-1)^{\ell}\right )+(k-1)^{\ell}\left ((k+1)^{\ell}-k^{\ell}\right )n\right )}{k^{\ell n}}.
\end{align*}
\end{theorem}
\begin{proof}
    The case $n=1$ is covered in Corollary \ref{crlComplete}. Let $n\geq 2$ and denoted by $S_n$ the sum of the number of blocks across all $k$-colored partition in $K_{\ell}\times P_n$. The claim is that
    \begin{align}
    \label{expFormula}
      k^{\ell n}\mathbb{E}[X_{\Til_2}(K_{\ell}\times P_n)]=k^{\ell n-(2\ell -1)}\left (\left (k^{2\ell}-(k^2-1)^{\ell}\right )+(k-1)^{\ell}\left ((k+1)^{\ell}-k^{\ell}\right )n\right).  
    \end{align}
    Let $S_n$ denote the right-hand side of  \eqref{expFormula}. We can express $S_n$  as
    \begin{align*}
        S_n &= k^{\ell}S_{n-1}+k^{\ell n-(2\ell -1)}(k-1)^{\ell}\left ((k+1)^{\ell}-k^{\ell}\right )\\
        &= k^{\ell}S_{n-1}+k\cdot k^{\ell (n-2)}(k-1)^{\ell}\left ((k+1)^{\ell}-k^{\ell}\right ).
    \end{align*}

Combinatorially, this decomposition separates two cases:  one where a block starts at the last slice and one where it does not. In the latter case, we obtain $k^{\ell}S_{n-1}$, as the number of possible colorings of the last slice is $k^{\ell}$, which means that each block is counted $k^{\ell}$ times. For the former case, where blocks are created in the new slice, we only need to consider the last two slices. The new block will appear $k^{\ell(n-2)}$ times, which accounts for all possible colorings of $K_{\ell}\times P_{n-2}$. Let $m$ denote the number of vertices that form this new block in the last slice. We can choose these vertices in $\binom{\ell}{m}$ ways and color them in one of $k$ colors.  The corresponding vertices on the previous slice must be colored with a different color,  giving  $(k-1)^{m}$ possibilities. The remaining vertices of the last slice have to be colored with a different color of the new block, giving $(k-1)^{\ell -m}$, and the corresponding vertices of the second to the last slice can be freely colored, producing $k^{\ell -m}$ possibilities. 

Summing over all  possible values of $m$, where $1\leq m\leq \ell$, we have $$k\cdot k^{\ell (n-2)}(k-1)^{m+(\ell-m)}\sum _{m = 1}^{\ell}\binom{\ell}{m}k^{\ell -m} = k\cdot k^{\ell (n-2)}(k-1)^{\ell}\left ((k+1)^{\ell}-k^{\ell}\right ).$$
By adding the contributions from both cases, we complete the induction step, proving the formula for $n$.
\end{proof}

\section{Colored Partition on the Star Graph}

\begin{figure}[ht!]
    \centering
    \includegraphics[width=0.75\linewidth]{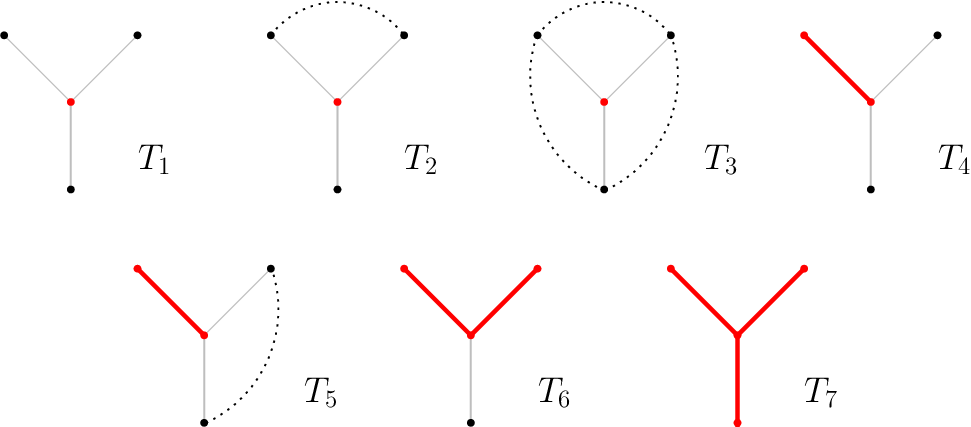}
    \caption{Each of the seven possible configurations for the  last slice of $K_{1,3}\times P_n$.}
    \label{figStarGraph}
\end{figure}

We can apply the technique used in the previous section to count the number of colored partitions for other graphs. Specifically, we enumerate the $2$ colored partitions for the product $K_{1,3}\times P_n$, where $K_{1,3}$, known as the \textit{star graph},  is defined as in Section \ref{sectionBipartite}. Figure \ref{figStarGraph} illustrates all possible configurations for the last slice of $K_{1,3}\times P_n$. In this figure, the dotted lines indicate that the two vertices they  connect belong to the same block, which is determined by a path that passes through other slices of the graph.

The relations  between the variables $T_i(x,y)$, which represent the bivariate generating functions for the $2$-colored partitions ending in the configuration $T_i$, can be encoded in the following $7\times 7$ matrix:
\begin{align*}
    \mathcal{R} = \left(\begin{array}{rrrrrrr}
y^{4} + 1 & 3 y^{4} & y^{4} & 3 y^{3} + 3 y & 3 y^{3} & 3 y^{2} & y^{3} \\
0 & 1 & 0 & 0 & y & y^{2} & 0 \\
0 & 0 & 1 & 0 & 0 & 0 & y \\
y^{2} + 1 & 3 y^{2} + 2 & y^{2} & y^{3} + 4 y + 1 & y^{3} + 4 y & 3 y^{2} + 2 y & y^{2} \\
0 & 1 & 1 & 0 & 1 & 1 & y \\
2 & y + 5 & y + 1 & \frac{3 y^{2} + 2 y + 1}{y} & 3 y + 3 & y^{2} + 2 y + 3 & 2 y \\
\frac{y^{2} + 1}{y^{2}} & \frac{3 y + 3}{y} & 2 & \frac{3 y + 3}{y} & 6 & 6 & y + 1
\end{array}\right).
\end{align*}

\begin{figure}[h!]
        \centering
        \includegraphics[width=0.6\linewidth]{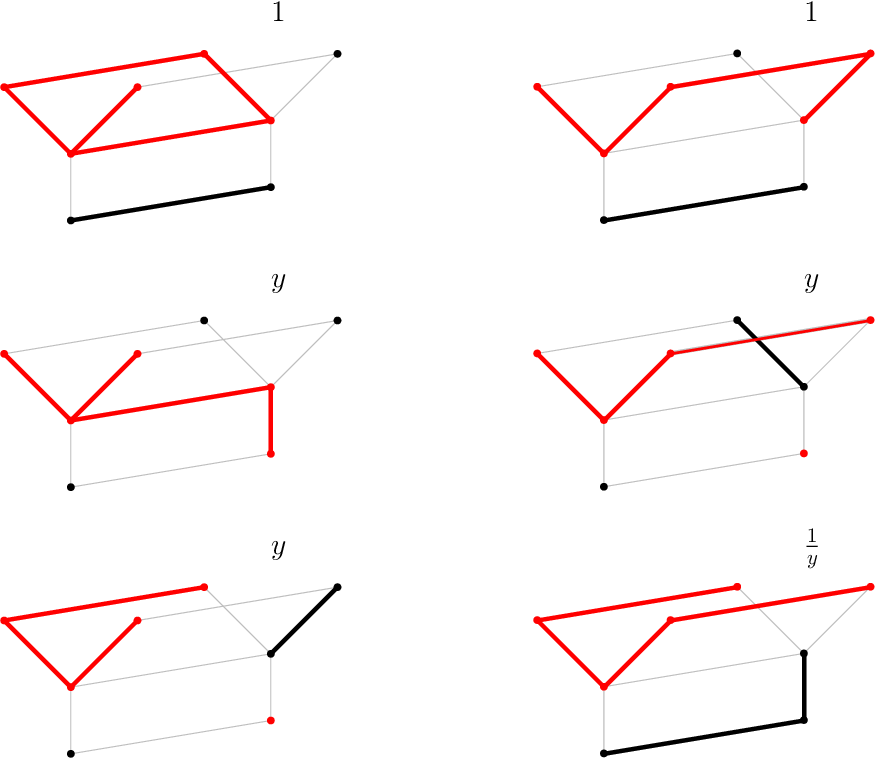}
        \caption{$T_6$ on top of all possible configurations $T_4$.}
        \label{fig:R64}
    \end{figure}

For example, consider the entry $\mathcal{R}_{6,4}$, which corresponds to adding the configuration $T_6$ to $K_{1,3}\times P_{n-1}$, where the last slice is one of six possible configurations $T_4$. Figure \ref{fig:R64} illustrates all six possible configurations. The generating function corresponding to these configurations is given by 
$$2+3y+\frac{1}{y} = \frac{3y^2+2y+1}{y}.$$

Let $\texttt{t} = (T_1(x,y),T_2(x,y),T_3(x,y),T_4(x,y),T_5(x,y),T_6(x,y),T_7(x,y))$. Using the matrix $\mathcal{R}$, the problem can be formulated as solving the system $\texttt{t} = \texttt{b}+x\mathcal{R}\texttt{t}$,
where $\texttt{b}=(xy^4,0,0,xy^3,0,xy^2,xy)$ is the vector that contains the base cases. To find $\texttt{t}$, we first compute $(I-x\mathcal{R})^{-1}$, where $I$ is the identity matrix. Multiplying this inverse matrix by   $\texttt{b}$  gives us the generating functions for each one of the last possible slices. Substituting these results into Equation \eqref{eqSliceStar} yields  the generating function:
\begin{align}
    \label{eqSliceStar}
    T(x,y) = 2\left (T_1(x,y)+3T_2(x,y)+T_3(x,y)+3T_4(x,y)+3T_5(x,y)+3T_6(x,y)+T_7(x,y)\right ).
\end{align}

\begin{theorem}
    The bivariate generating function for the 
$2$-colored partitions of the product graph $K_{1,3} \times P_n$ is given by $T(x,y) =p(x,y)/q(x,y)$,
    where  \footnotesize
    \begin{align*}
        p(x,y) &=(2y^4 + 6y^3 + 6y^2 + 2y)x + (-8y^7 + 4y^6 - 28y^5 + 40y^4 - 44y^3 + 4y^2 - 16y)x^2+\\
        & (14y^9 - 22y^8 - 4y^7 + 104y^6 - 106y^5 - 32y^4 + 90y^3 - 52y^2 + 40y)x^3+\\
        & (-8y^{10} - 8y^9 + 124y^8 - 204y^7 + 52y^6 + 64y^5 + 16y^4 - 80y^3 + 92y^2 - 48y)x^4 +\\
        &(-10y^{11} + 4y^{10} + 110y^9 - 268y^8 + 112y^7 + 384y^6 - 662y^5 + 468y^4 - 120y^3 - 48y^2 + 30y)x^5 +\\
        &(16y^{11} - 28y^{10} - 156y^9 + 656y^8 - 1028y^7 + 708y^6 - 316y^4 + 168y^3 - 12y^2 - 8y)x^6 +\\
        &(8y^{11} - 58y^{10} + 228y^9 - 562y^8 + 852y^7 - 756y^6 + 340y^5 - 26y^4 - 36y^3 + 10y^2)x^7,
    \end{align*} \normalsize
    and \footnotesize
    \begin{align*}
        q &= 1 + (-y^4 - y^3 - y^2 - 7y - 9)x + (y^7 + y^6 + 5y^5 + 11y^4 - 4y^3 - 6y^2 + 14y + 28)x^2 +\\ &(-y^9 - 3y^8 + y^7 + 5y^6 - 22y^5 - 11y^4 + 7y^3 + 54y^2 - 18y - 44)x^3 +\\
        &(7y^9 - 17y^8 + 2y^7 + 20y^6 - 32y^5 + 45y^4 + 42y^3 - 105y^2 - y + 39)x^4 +\\
        &(y^{11} + 4y^{10} - 24y^9 + 47y^8 - 28y^7 - 62y^6 + 167y^5 - 125y^4 - 50y^3 + 83y^2 + 6y - 19)x^5 +\\
        &(y^{11} - 24y^{10} + 94y^9 - 122y^8 - 61y^7 + 365y^6 - 409y^5 + 116y^4 + 116y^3 - 91y^2 + 11y + 4)x^6 +\\ &(-3y^{11} + 23y^{10} - 74y^9 + 95y^8 + 45y^7 - 289y^6 + 355y^5 - 183y^4 + 18y^3 + 18y^2 - 5y)x^7.
    \end{align*}
\end{theorem}
\normalsize
\begin{corollary}
   The expected number of blocks when uniformly and independently coloring the graph $K_{1,3}\times P_n$ with $3n$ vertices using two colors is given by
    $$\mathbb{E}\left [X_{\Til_2}(K_{1,3}\times P_n)\right ]=\frac{2254219}{1411200}+\frac{6}{49}2^{-3n}-\frac{2}{225}2^{-4n}+\frac{11933}{13440}n.$$
\end{corollary}

To count the size of the system of equations for a general enumeration of  $2$-colored partitions for the graph $K_{1,m}\times P_n$ based on the last possible slices, we use the function $p(n)$ that counts the number of integer partitions of $n$. For a given $m$,  the number of different last slices corresponds to the sum of the number of partitions of integers up to $m$. Specifically,  the number of equations required to solve the system is given by
$$ \sum _{\ell = 0}^{m}p(\ell).$$ 
The partition function comes from all possible ways to connect points colored the opposite way of the color given to the vertex with degree $m$.
For example, for $m = 3$, we have $p(0)+p(1)+p(2)+p(3) = 1+1+2+3= 7$ possible last slices shown in Figure \ref{figStarGraph}.


\begin{thebibliography}{20}


\bibitem{belcastro} s-m.~belcastro. Domino tilings of $2\times n$ grids (or perfect matchings of grid graphs) on surfaces. J.  Integer Seq. \textbf{26} (2023), Article 23.5.6.

\bibitem{Bodini} O.~Bodini. On the strange kinetic aesthetic of rectangular shape partitions. Pure Math. Appl. (PU.M.A.) \textbf{30} (2022), 37--44.

\bibitem{DP} T.~Do\v{s}lic and L.~Podrug. Sweet division problems: from chocolate bars to honeycomb strips and back.  Accepted in the American Mathematical Monthly, (2024).
 

\bibitem{Kas}
P.~W.~Kasteleyn. Dimer statistics and phase transitions. J. Mathematical  Phys. {\bf 4} (1963), 287--293.
  
\bibitem{Mansour}
R.~Mansour. Counting clusters in a coloring grid. Discrete Math. Lett. \textbf{5} (2021), 20--23.

 \bibitem{RamVil} J.~L.~Ram\'irez and D.~Villamizar. Colored random tilings on grids. J. Autom. Lang. Comb. (2024).

 \bibitem{RamVil2} J.~L.~Ram\'irez and D.~ Villamizar. Counting colored tilings on grids and graphs. In:  Proceedings of the 13th edition of the conference on Random Generation of Combinatorial Structures. Polyominoes and Tilings (GASCom 2024), Bordeaux, France, 24-28th June 2024. Electron. Proc. Theor. Comput. Sci. (EPTCS) \textbf{403} (2024), 164--168.

\bibitem{Richey} J.~Richey. Counting clusters on a grid. Undergraduate Honors Thesis. Dartmouth College, 2014.


\bibitem{RolUg} N.~Rolin and A.~Ugolnikova. Tilings by $1\times 1$ and $2\times 2$.  RAIRO-Theor. Inf. Appl. \textbf{50} (2016), 105--116


\bibitem{Stanley}
R.~P.~Stanley. \emph{Algebraic Combinatorics: Walks, Trees, Tableaux, and More}. Springer, 2013.

\bibitem{TM}
H.~N.~V.~Temperley and M.~E.~Fisher. Dimer problem in statistical
  mechanics---an exact result.  Philos. Mag. {\bf 6} (1961),
  1061--1063.
  
\end{thebibliography}
\end{document}